\documentclass{amsart}

\newtheorem{rem}{Remark}[section]
\newtheorem{teo}{Theorem}[section]
\newtheorem{defin}{Definition}[section]
\newtheorem{prop}{Proposition}[section]
\newtheorem{coro}{Corollary}[section]
\newtheorem{lemma}{Lemma}[section]

\def\proof{{\bf Proof:}\ }
\def\endproof{{\mbox{}\nolinebreak\hfill\rule{2mm}{2mm}\par\medbreak} }

\def\eq#1{(\ref{#1})}
\def\neweq#1{\begin{equation}\label{#1}}
\def\endeq{\end{equation}}

\def\bbbr{\mathbb{R}}
\def\bbbn{\mathbb{N}}
\def\half{\frac12}
\def\supp{{\rm supp}}
\def\oint{\int_\Omega}
\def\a{\alpha}

\def\eps{\varepsilon}
\def\l{\lambda}

\def\O{\Omega}
\def\g{\gamma}

\def\o{\omega}
\def\spz{H^1_0(\O)}
\def\spk{(H^1_0(\O))^k}
\def\N{{\mathcal N}}

\def\dis{\displaystyle}

\title{On a class of optimal partition problems  \\
related to the Fu\v{c}\'\i k spectrum and to\\ the monotonicity
formulae }

\author{Monica Conti, Susanna Terracini and Gianmaria Verzini}
\address{Monica Conti and Gianmaria Verzini
\newline\indent
Politecnico di Milano
\newline\indent
Dipartimento di Matematica
\newline\indent
Via Bonardi 9, I-20133 Milano, Italy
\newline
\newline\indent
Susanna Terracini
\newline\indent
Universit\`a  di Milano Bicocca
\newline\indent
Dipartimento di Matematica ed Applicazioni
\newline\indent
Piazza dell'Ateneo
Nuovo 1, 20126 Milano, Italy}
\date{\today}
\thanks{Work partially supported by MIUR
project ``Metodi Variazionali ed Equazioni Differenziali Non
Lineari''}

\begin{document}

\maketitle \pagestyle{plain}

\begin{abstract}
In this paper we give an unified approach to some questions
arising in different fields of nonlinear analysis, namely: (a) the
study of the structure of the Fu\v{c}\'\i k spectrum and (b)
possible variants and extensions of the monotonicity formula by
Alt--Caffarelli--Friedman \cite{acf}. In the first part of the
paper we present a class of optimal partition problems involving
the first eigenvalue of the Laplace operator. Beside establishing
the existence of the optimal partition, we develop a theory for
the extremality conditions and the regularity of minimizers. As a
first application of this approach, we  give a new variational
characterization of the first curve of the Fu\v{c}\'\i k spectrum
for the Laplacian, promptly adapted to more general operators. In
the second part we prove a monotonicity formula in the case of
many subharmonic components and we give an extension to solutions
of a class of reaction--diffusion equation, providing some
Liouville--type theorems.
\end{abstract}


\section{Introduction and statement of the results}

Let $\Omega\subset\bbbr^N$ be a connected, open bounded domain
with regular boundary $\partial\Omega$. For any open $\omega\subset
\Omega$, let $\lambda_1(\omega)$ denote the first eigenvalue of
the Laplace operator in $H^1_0(\omega)$, namely
$$\lambda_1(\omega)=\min_{u\in H^1_0(\omega)\atop u\not\equiv 0} \frac{\int_\omega |\nabla
u(x)|^2dx}{\int_\omega|u(x)|^2dx}.$$
For a fixed $p>0$, let us consider the following
class of optimal partition problems
\neweq{model}
\inf_{{\mathcal P}_k}\frac{1}{k}\sum_{i=1}^k
\big(\lambda_1(\omega_i)\big)^p,
\endeq
where the minimization is taken over the class of partitions in
$k$ disjoint, connected, open subsets of $\O$
\neweq{partizioni}
{\mathcal P}_k:=\left\{(\omega_1,\dots,\omega_k)\subset\O:
\;\omega_i\mbox{ is open and
connected},\;\omega_i\cap\omega_j=\emptyset\;{\rm if}\;i\neq
j\right\}\;.
\endeq In this paper we shall investigate different aspects of
problem \eqref{model} and of the analogous optimal partition
problem on the spheres of $\bbbr^N$.

The motivation of our interest in problems of this type is that
they are beyond different and relevant questions of nonlinear
analysis. For instance, the proof of the monotonicity formula by
\cite{acf}, relies on the determination of the value
$$\inf_{\omega\subset S^{N-1}}
\g(\lambda_1(\omega))+\g(\lambda_1(S^{N-1}\setminus\omega)),$$
where $\g$ denotes a suitable function (namely,
$\g(s)=\sqrt{(N/2-1)^2+s}$), and $S^{N-1}$ denotes the boundary of
the unit ball in $\bbbr^N$. Similar optimal partition problems
arise in connection with the phenomenon of the spatial segregation
in reaction--diffusion systems, as shown in \cite{ctv3,ctv4}.

In the present paper, we develop a theory for \eq{model} and we
give applications in two independent directions. First, by
studying problem \eqref{model} for partitions of connected domains
of $\bbbr^N$, we obtain  a new variational characterization of the
first curve of the Fu\v{c}\'\i k spectrum of the Laplacian, that
can be promptly adapted to the case of the $p$--Laplacian and even
to more general notions of spectrum. On the other hand, by
exploiting partitions of $S^{N-1}$, we shall prove some
monotonicity formulae related to \cite{acf}, for the case of many
components.\\ A crucial tool will be the theory already developed by the authors in \cite{ctv2},
concerning the regularity of solutions and the properties of the free boundary,
in connection  with certain classes of
optimal partition problems involving nonlinear eigenvalues. More
recently, in \cite{ctv3}, the theory is shown to apply in the more
general context of $k$--uple of functions with mutually disjoint
supports and belonging to functional classes characterized by
suitable differential inequalities.

Let us now describe the structure of the paper and our main
results with some details.

The first part of the paper is devoted to the study of optimal
partition problems of the general type of \eq{model}. A different
class of optimal partition problems, with an area constraint, has
been studied in \cite{bu1,bu2}. To begin with, we seek the best
partition of $\Omega$ with respect to the infimum defined in
\eqref{model}. Then, our main goal is to derive the extremality
sufficient conditions and to prove some qualitative properties
both of the optimal sets $\omega_i$'s and of the eigenfunctions
associated to $\lambda_1(\omega_i)$. Our first result in this
direction is:

\begin{teo}\label{extr_intr}
There exists $(\omega_1,\dots,\omega_k)\in {\mathcal P}_k$
achieving \eq{model}. Furthermore, if $\phi_1,\dots,\phi_k$ are
the associated eigenfunctions normalized in $L^2$, then, there
exist $a_i\in\bbbr$ such that the functions $u_i=a_i\phi_i$
verifies in $\Omega$ the differential inequalities
\begin{itemize}
\item[1.] $-\Delta u_i\leq \lambda_1(\omega_i)u_i$,
\item[2.] $-\Delta\left( u_i-\sum_{j\neq i}u_j\right)\geq \lambda_1(\omega_i)u_i-\sum_{j\neq i} \lambda_1(\omega_j)u_j$.
\end{itemize}
\end{teo}
As a consequence, the main results of \cite{ctv3} apply, providing
the regularity of the minimizing   $k$--uple $U=(u_i)_i$ and some
qualitative features of the interfaces $\partial
\omega_i\cap\partial \omega_j$, together with an asymptotic
expansion of $U$ at multiple intersection points in dimension $2$.
Finally, we extend all the previous considerations to the limiting
case in \eq{model} (as $p\to\infty$):
\neweq{fucik_intr}
\inf_{(\omega_i)\in{\mathcal P}_k}\max_{i=1,\dots,k}
\{a_i\lambda_1(\omega_i)\},
\endeq
where $a_i$ are given positive weights. This last problem  is
connected to  the Fu\v{c}\'\i k spectrum $\mathcal F$ of the
Laplace operator with Dirichlet boundary conditions, which was
defined in \cite{d0,f} as the set of pairs $(\lambda,\mu)$ such
that the problem
$$
-\Delta u=\lambda u^+-\mu u^-\;,\qquad u\in H^1_0(\Omega),
$$
has a non zero solution (here we use the standard notation $u^\pm:=
\max\{\pm u,0\}$).

In fact, to each solution of \eq{fucik_intr} there corresponds an
element of $\mathcal F$, provided either $k=2$ or the boundary of
the supports $\partial\omega_i$ do not have multiple intersection
points.  In this way we can give a new variational
characterization of the first curve of the Fu\v{c}\'\i k spectrum,
in terms of an optimal partition of eigenvalues. More precisely,
let
$$
c(r):=\inf_{(\omega_i)\in{\mathcal P}_2}
\max\{r\lambda_1(\omega_1),\lambda_1(\omega_2)\}.
$$
Then we have
\begin{teo}\label{fuc_intr}
For all $r>0$, there exists $u\in H^1_0(\Omega)$ such that
$(\{u^+>0\},\{u^->0\})$ achieves $c(r)$. Furthermore, the pair
$$(\lambda_1(\{u^+>0\}),\lambda_1(\{u^->0\}))=(r^{-1}c(r),
c(r))$$
 belongs to the Fu\v{c}\'\i k spectrum and it represents
  the first nontrivial intersection between ${\mathcal F}$ and the line of slope $r$.
\end{teo}
Variational formulations of the first curve have been given as
min--max or constrained minimum, starting form the
one--dimensional periodic problem in \cite{dfr}; the general case
$N\geq 1$ has been first studied in the paper of De Figuereido and
Gossez \cite{dfg}. More recently, new characterizations were
proposed in \cite{cdfg} and \cite{hr}, covering the  Fu\v{c}\'\i k
spectrum of the $p$--laplacian;  theoretical and numerical studies
have been carried out in \cite{br1,br2,hr}. In our opinion,  our
characterization is of interest from both the theoretical and the
computational points of view. Indeed, it admits straightforward
extensions to many other nonlinear operators, such as the
$p$--laplacian and it could be easily modified in order to apply
to general boundary conditions; moreover it can be easily
implemented numerically using, for instance, a steepest descent
method. A formulation related to our Theorem \ref{fuc_intr} was
given in \cite{hr}, but could not cover the full first curve of
the spectrum.

Finally, as far as the dimension 1 is concerned, we can
characterize an infinity of curves belonging to the Fu\v{c}\'\i k
spectrum, in the following way:
\begin{teo}\label{altre-curve}
Let $N=1$, $k\geq1$ and define, for all $r>0$
$$
c_{k+1}(r)=\inf_{a=t_0<t_1<...<t_k<t_{k+1}=b}
\max_{i}\{r\lambda_1(t_{2i+1}-t_{2i}),\lambda_1(t_{2i+2}-t_{2i+1})\}.$$
Then the pair
$$(r^{-1}c_{k+1}(r),c_{k+1}(r))$$
belongs to ${\mathcal F}$.
\end{teo}
The second part of the paper is devoted to the study of the
monotonicity formulae. The {\it monotonicity lemma}  was
originally stated  by Alt, Caffarelli and Friedman in \cite{acf}
in the following way:
\begin{lemma}[The monotonicity formula]\label{clax_intr}
Let $(w_1,w_2)\in (H^1(\Omega))^2$ be non negative, continuous,
subharmonic functions in a ball $B(x_0,\bar r)\subset\Omega$ (i.e.
$-\Delta w_i\leq 0$ in distributional sense). Assume that
$w_1(x)w_2(x)=0$. Assume that $x_0\in\partial(\supp(w_i))$ for
$i=1,2$. Define
$$\Phi(r) =\prod_{i=1}^{2}\frac{1}{r^2}\int_{B(x_0,r)}\frac{|\nabla
w_i(x)|^2}{|x-x_0|^{N-2}}dx.$$ Then $\Phi$ is a non decreasing
function in $[0,\bar{r}]$.
\end{lemma}

Since its very first publication, the monotonicity formula was
shown to be a powerful tool in proving many local results in the
theory of free boundaries. Our first objective consists in
developing a variant of Lemma \ref{clax_intr} for the case of many
subharmonic densities having mutually disjoint supports. To this
aim we consider the optimal partition value
$$
\beta(k,N):= \inf_{{\mathcal P}(k,N)}\frac{2}{k}\sum_{i=1}^k \sqrt{\lambda_1(\omega_i)},
$$
where the minimization is taken over all possible partitions in
$k$ disjoint parts of the unit sphere $S^{N-1}$. We shall prove:
\begin{lemma}\label{monot.generale_intr}
Let $\Omega\subset \bbbr^N, N\geq 2$. Let $w_1,\dots,w_h\in
H^1(\Omega)$ be non negative subharmonic functions in a ball
$B(x_0,\bar r)\subset\Omega$ (i.e. $-\Delta w_i\leq 0$ in
distributional sense). Assume that $w_i(x)w_j(x)=0$ a.e. if $i\neq
j$ and that $x_0\in\partial(\supp(w_j))\cap\Omega$ for all
$j=1,...,h$. Define
$$\Phi(r) =\prod_{i=1}^{h}\frac{1}{r^{\beta(h,N)}}
\int_{B(x_0,r)}|\nabla w_i(x)|^2dx.$$
Then $\Phi$ is a non decreasing function in $[0,\bar{r}]$.
\end{lemma}
In the recent years, several papers have shown the existence of a
strong connection between some free boundary problems and the
spatially segregated limits of competition--diffusion systems, as
the interaction rates tend to infinity. This asymptotic study has been
carried out in \cite{d1,dd3,dg1,dhmp,kl,lm}. The link with the optimal
partitions was examined by the authors in \cite{ctv3,ctv4}. This
motivates the interest of extending the monotonicity lemma to the
case of an arbitrary number of densities whose supports need not
to be mutually disjoint, but, instead, satisfy a system of
competition--diffusion equations. In particular, as a prototype we
consider the following system:
\neweq{k_nonsimm_intr}
\left\{\begin{array}{rcll}
 -\Delta u_i(x)&=&\dis -u_i(x)\sum_{j\neq i}a_{ij}u_j(x),\quad   &x\in\bbbr^N,\\
  u_i(x)&\geq&0,                                           &x\in\bbbr^N,
  \end{array}\right.
 \endeq
for all $i=1,\dots,k$, where $a_{ij}>0$.
Then we shall prove the following monotonicity result:
\begin{lemma}\label{monot.formula.N_intr}
Let $N\geq 2$ and let $(u_1,\dots,u_k)$ be a solution of \eq{k_nonsimm_intr}
such that $u_i>0$ for all $i$. Let
$h\leq k$ be any integer, let $h'<\beta(h,N)$ and define
$$
\Phi(r)=\prod_{i=1}^{h}\frac{1}{r^{h'}}\int_{B(0,r)}\left(|\nabla
u_i(x)|^2+u_i^2(x)\sum_{j\neq i}a_{ij}u_j(x)\right)dx\;.
$$
Then there exists $r'=r(h')$ such that $\Phi$ is an increasing
function in $[r',\infty)$.
\end{lemma}

The above perturbed monotonicity formula will turn out to be the
key point in proving an a priori growth estimate on the non
trivial solutions of  \eq{k_nonsimm_intr}. The subsequent
Liouville type result will be exploited, in a forthcoming paper
\cite{ctv4}, in order to prove the equi--h\"olderianity of
solutions of competition--diffusion systems, when the
interspecific competition rate tends to infinity.

\section{Optimal partition problems involving linear
eigenvalues}

Let $p$ be a positive real number; let $k\geq 2$ be a fixed
integer. In this section we consider the problem of finding a
partition of $\Omega$  that achieves
\neweq{problema}
\inf_{(\omega_i)\in{\mathcal P}_k}\sum_{i=1}^k
(\lambda_1(\omega_i))^p,
\endeq
where ${\mathcal P}_k$ is the set of all the possible partitions
of $\O$ in $k$ connected, open subsets (see \eq{partizioni}). Let
us recall that $\lambda_1(\omega)$ denotes the first eigenvalue of
$-\Delta$ in $H^1_0(\o)$; the associated eigenspace is one
dimensional (if $\omega$ is connected), and the associated
eigenfunction does not change its sign.

As a first step of our investigation, we relax problem
\eq{problema} in the following way. For any measurable
$\omega\subset\Omega$, let $\lambda_1(\omega)$ denote
$$
\lambda_1(\omega):=\inf\left\{\frac{\int_\omega |\nabla
u(x)|^2dx}{\int_\omega|u(x)|^2dx}:\,u\in H^1_0(\Omega),\,u=0
\mbox{ a.e. on }\O\setminus\o,\,u\not\equiv 0\right\}.
$$
When $\o$ is open, then $\l_1(\o)$ is the classical first
eigenvalue of the Laplace operator in $H^1_0(\o)$. For this
reason, with a slight abuse of notation, we shall name, for any
arbitrary measurable set $\o$, $H^1_0(\o):=\{u\in
H^1_0(\Omega),\,u=0 \mbox{ a.e. on }\O\setminus\o\}$
(incidentally, we observe that possibly $H^1_0(\o)\equiv\{0\}$,
and consequently we define in a standard way $\l_1(\o)=+\infty$).
Following this line, let us introduce the class of relaxed
partitions
$$
{\mathcal P}^*_k:=\left\{(\omega_1,\dots,\omega_k)\subset\O:
\;\omega_i\mbox{ is
measurable},\;\omega_i\cap\omega_j=\emptyset\;{\rm if}\;i\neq
j\right\}\;
$$
and the relaxed minimization problem
\neweq{problema*}
\inf_{(\omega_i)\in{\mathcal P}^*_k}\sum_{i=1}^k
(\lambda_1(\omega_i))^p,
\endeq
Clearly the value in \eq{problema*} is smaller or equal to the one
in \eq{problema}. Nevertheless, at the end of this section the
reader will see that \eq{problema} and \eq{problema*} are in fact
equivalent.

By taking into account the variational characterization of
$\lambda_1$, the infimum value in \eq{problema*} becomes
\neweq{probweak}
 \inf_{u_i\in H^1_0(\O)\setminus\{0\}\atop{u_i\cdot u_j=0}}\sum_{i=1}^k
\left(\frac{\int_\Omega |\nabla u_i(x)|^2dx}{\int_\Omega|u_i(x)|^2dx}\right)^p.
\endeq
We have

\begin{rem}\label{deboleesistenza}
There exists a $k$--uple $(\phi_1,\dots,\phi_k)$, with
$\|\phi_i\|_2=1$, that achieves the value \eq{probweak}. Moreover,
if $\omega_i=\{x\in\O:\phi_i(x)>0\}$, then $\sum
(\lambda_1(\omega_i))^p$ achieves \eq{problema*}. {\rm Indeed, we
can minimize the functional
$$
{\mathcal E}(u_1,...,u_k)=\sum_{i=1}^k
\left(\frac{\int_\Omega |\nabla
u_i(x)|^2dx}{\int_\Omega|u_i(x)|^2dx}\right)^p,
$$
among $k$--uples of $H^1_0$ functions, subject to the constraint
that $u_i\geq 0$ for all $i$ and $u_i\cdot u_j=0$ if $i\neq j$.
Since the above functional is weakly lower semicontinuous, and the
constraint is locally weakly compact, the direct method of the
calculus of variations applies. This immediately provides the
existence of a $k$--uple of functions $(\phi_i)$ (that we can
obviously assume normalized in $L^2$), having disjoint supports,
which achieves \eq{probweak}. As a consequence, letting
$\omega_i=\{\phi_i>0\}$, we find a solution of \eq{problema*}.}
\end{rem}

The following part of the section (Subsections \ref{esistenza},
\ref{subsup}) is devoted to the study of the properties of the
minimizers of \eq{problema*} in the case $p\neq1$ (the case $p=1$
will be considered in Section \ref{p=1} through a limiting
procedure). Finally, in Subsection \ref{ultimasubsection} we will
show the equivalence between \eq{problema} and \eq{problema*}.

\subsection{An auxiliary variational problem (for $p\neq1$)}\label{esistenza}

In this section we prove the equivalence of  our minimal partition
problem \eq{problema*} with a min--max value for a certain
functional defined on $(H^1_0)^k$. Let us start the description of
the appropriate variational setting with some definitions:
 first, let $p\neq 1$ be fixed and let $q$ be the dual exponent
of $p$, $q=\frac{p}{p-1}$. Note that $q\in (-\infty,0)\cup
(1,\infty)$.
For $u\in\spz$ and $U:=(u_1,\dots,u_k)\in\spk$ we
define the functionals
$$
\begin{array}{rcl} \dis J^*(u)&:=&\frac{1}{2}\int_\O |\nabla
u(x)|^2dx - \frac{1}{2q}\left(\int_\O u(x)^2dx\right)^q\\\ \dis
J(U)&:=&\sum_{i=1}^k J^*(u_i).\end{array}
$$
We are interested in studying $J$ restricted to the $k$--uples of
functions where the different components have disjoint support,
hence we define
$$
{\mathcal H}=\{U:=(u_1,\dots,u_k)\in\spk\;:\; u_i\cdot u_j=0,{\rm  a.e. on }\,\O\,{\rm for }\,i\neq
j\}.
$$
Note that $U\in \mathcal H$ implies
$J(U)=\sum_{i=1}^kJ^*(u_i)=J^*(\sum_{i=1}^k u_i).$

Next we define the Nehari manifolds associated to $J^*$ and $J$
$$\begin{array}{rcl}
\N(J^*)&:=&\{u\in\spz:\,u\geq0,\,u\not\equiv0,\,\nabla J^*(u)\cdot
u=0\}\\
\N_0&:=&(\N(J^*))^k \cap {\mathcal H},
\end{array}$$
and the value
\neweq{ck}\begin{array}{c}
c_q:=\dis\inf\left\{ J(U):\;U\in\N_0\right\}.
 \end{array}\endeq

We shall prove the following
\begin{teo}\label{optimal}
Let either $q<0$ or $q>1$, and let $\frac{1}{p}+\frac{1}{q}=1$.
Then
\neweq{pb_debole}c_q\equiv\inf_{(\omega_i)\in{\mathcal P}^*_k}
\frac{q-1}{2q}\sum_{i=1}^k \lambda_1(\omega_i)^p.\endeq Moreover,
there exists a $k$--uple of functions $U:=(u_1,\dots,u_k)\in\N_0$
such that $U$ achieve \eq{ck} and each $u_i$ solves
\neweq{equazioni}-\Delta u_i(x)=\left(\int_\O
u_i^2\right)^{q-1}u_i(x) \qquad\mbox{ when } u_i(x)>0.\endeq In
particular, for all $i$, there holds
\neweq{lambdaq}\lambda_1(\{u_i>0\})\equiv \left(\int_\O
u_i^2\right)^{q-1}.\endeq
\end{teo}
\proof
first, let $u\neq 0$,
$\omega:=\{u>0\}$ and consider $J^*$ restricted to the line
$t\longmapsto tu$, namely $g(t)=J^*(tu)$. It turns out that, for
$t\neq 0$, $g'(t_u)=0$ iff
$${ t_u}^{2q-2}=\frac{\int_\omega |\nabla
u(x)|^2dx}{\int_\omega|u(x)|^2dx},$$ and thus
$$
J^*(t_u u)=\frac{q-1}{2q}\left(\frac{\int_\omega |\nabla
u(x)|^2dx}{\int_\omega|u(x)|^2dx}\right)^\frac{q}{q-1}
$$
As a consequence
\neweq{bart}
\inf_{u\in H^1_0(\omega)\atop u\not\equiv 0}J^*(\bar t_u u)
=\frac{1}{2p}\lambda_1(\omega)^p.
\endeq
Now some differences are induced by the value of $q$. In the case
$q>1$, the function $g$ has a local minimum at the origin and
$\lim_{t\to\infty}g(t)=-\infty$, hence $ \bar{t}_u$ is a local maximum.
As a consequence the value $c_q$ has the equivalent
characterization
$$
c_q=\inf_{U\in{\mathcal H}}\sum_{i=1}^k J(\bar{t}_{u_i}u_i)=
\inf_{U\in{\mathcal H}}\max_{\{(t_1,\dots,t_k):t_i>0\}}
\sum_{i=1}^k J(t_iu_i),\qquad\text{if }\;q>1.
$$
On the other side, if $q<0$ the corresponding function $g(t)$ is
such that $\lim_{t\to\infty}g(t)=\lim_{t\to 0}g(t)=\infty$  and
hence $c_q$ is the global infimum of $J$
$$
c_q=\dis\inf_{U\in \mathcal H} J(U),\qquad\text{if }\;q<0.
$$
Taking into account these characterizations  for the value $c_q$
and \eqref{bart}, we immediately have \eqref{pb_debole}.

Now that we have proved the equivalence between the value $c_q$
and the value \eqref{problema*}, the existence of a minimizer $U$
for \eq{ck} follows from the existence of a minimal partition for
problem \eq{problema*}, as we discussed in Remark
\ref{deboleesistenza}. Then, standard critical point techniques
prove that $U$ is a critical point for $J$, i.e. $\nabla
J(U)\equiv 0$ in distributional sense. By computing this means
$$\int_\O \nabla u_i\nabla v+\left(\int_\O u_i^2\right)^{q-1}\int_\O u_iv=0
\qquad\forall v\in H^1_0(\O)$$
and proves \eqref{equazioni}.
Finally, \eqref{lambdaq} is obtained in light of \eqref{bart}.
\endproof

\subsection{The extremality conditions}\label{subsup}
In this section we prove
that the extremality condition stated in Theorem \ref{extr_intr}, are verified by the solutions of the auxiliary problem \eq{ck}. For easier notation we define, for every $i$,
$$
\widehat u_i:=u_i-\sum_{j\neq i}u_j.
$$

\begin{lemma}\label{super} Let $U=(u_1,\dots,u_k)$ be as in Theorem
\ref{optimal} and  $(\omega_1,...,\omega_k)$ be the corresponding
supports. Then the following differential inequalities hold in
$\O$
\begin{itemize}
\item[1.] $-\Delta u_i\leq \lambda_1(\omega_i)u_i$,
\item[2.] $-\Delta\widehat u_i\geq \lambda_1(\omega_i)u_i-\sum_{j\neq i} \lambda_1(\omega_j)u_j$.
\end{itemize}
\end{lemma}

\proof the argument is different according to the case that $q>1$ or
$q<0$ and it mimics the proof in  \cite{ctv2} for the case of nonlinear eigenvalues.
 For the
reader's convenience we report the proofs adapted to the actual
setting. \\
{\it The case $q<0$}.\\
Let us prove 1. We argue by contradiction, assuming the existence of an index $j$ such that the claim does not hold;
 that is,  there exists $0\leq\phi\in C^\infty_c(\Omega)$
such that
\neweq{aa}
\int_\Omega\big[\nabla
u_j\nabla\phi-\lambda_1(\omega_j)u_j\phi\big]dx>0.
\endeq
For $t>0$ very small we define a new test function
$V=(v_1,\dots,v_k)$, belonging to ${\mathcal H}$, as follows: $$
v_{i}=\left\{\begin{array}{ll}u_{i} &\mbox{ if $i\neq j$},
\\(u_i-t\phi)^+&\mbox{ if $i=j$.}\end{array}\right.$$
We claim that $V$ lowers the value of the functional $J$. We introduce
$G(s)=\frac{1}{2q}s^q$ and compute as follows
$$
\begin{array}{rcl}
J(V)-J(U)
&=&\int_{\Omega}\kern-3pt{\frac12}\big(|\nabla
(u_j-t\phi)^+|^2-|\nabla u_j|^2\big)
-G\big(\int_{\Omega}\kern-3pt ((u_j-t\phi)^+)^2\big)+
G\big(\int_{\Omega}\kern-3pt (u_j)^2\big)
 \\
&\leq&\kern-3pt\int_{\Omega}\kern-3pt\frac12\big(|\nabla
(u_j-t\phi)|^2-|\nabla u_j|^2\big)
+2tG'\big(\int_{\Omega}\kern-3pt
(u_j)^2\big)\int_{\Omega}\kern-3pt u_j\phi +o(t) \\ &\leq&
-t\int_\Omega\big[\nabla u_j\nabla\phi-
2G'\big(\int_{\Omega}\kern-3pt (u_j)^2\big) u_j\phi\big]
+o(t).
\end{array}
$$
Note that the last expression, when $t$ is sufficiently small, is
negative by \eq{aa}, since $\lambda_1(\omega_j)\equiv
2G'\left(\int_{\Omega}\kern-3pt (u_j)^2\right)$; hence, choosing
$t$ sufficiently small, we obtain the contradiction
$$J(V)-J(U)<0.$$
In order to prove 2., let $j$ and $0<\phi\in C^\infty_c(\Omega)$ such that
$$
\int_\Omega\Big[\nabla \widehat u_j\nabla\phi
-\Big(\lambda_1(\omega_j) u_j-\sum_{i\neq j}
\lambda_1(\omega_i)u_i\Big)\phi\Big]\,dx<0.
$$
Again, we show that the value of the functional can be lessen by
replacing $U$ with an appropriate new test function $V$. To this
aim we consider the positive and negative parts of $\widehat
u_j+t\phi$ and we notice that, obviously,
$$\{(\widehat u_j+t\phi)^->0\}\subset\{(\widehat u_j)^->0\}=\cup_{i\neq j}\{u_i>0\}\;.$$
Let us define $V=(v_1,\dots,v_k)$ in the following way:
$$ v_i=\left\{
\begin{array}{ll}
\left(\widehat u_j+t\phi\right)^+, &\mbox{if $i=j$},\\
\left(\widehat u_j+t\phi\right)^-\chi_{\{u_i>0\}}, &\mbox{if $i\neq j$}.
\end{array}\right.$$
Here and below $\chi_A$ denotes the characteristic function of the set $A$.
We compute as follows
\neweq{co7}
\begin{array}{rcl}
J(V)-J(U)&=&\sum_{i=1}^k\int_\Omega\frac12 \big(|\nabla
v_i|^2-|\nabla u_i|^2\big)dx- G\big(\int_{\Omega}\kern-3pt
(v_i)^2\big)+G\big(\int_{\Omega}\kern-3pt (u_i)^2\big)
\\
&=&\int_{\Omega}\frac12\big(|\nabla \widehat u_j+t\phi|^2-|\nabla
\widehat u_j|^2\big)dx -G\big(\int_{\Omega}\kern-3pt ((\widehat
u_j+t\phi)^+)^2\big)+ G\big(\int_{\Omega}\kern-3pt
(u_j)^2\big)-
\\&-&\sum_{i\neq j}\big[
G\big(\int_{\{u_i>0\}}\kern-3pt ((\widehat
u_j+t\phi)^-)^2\big)-G\big(\int_{\Omega}\kern-3pt
(u_i)^2\big)\big]
\\&=&
t\int_{\Omega}\nabla \widehat u_j\nabla\phi-
2tG'\big(\int_{\Omega}\kern-3pt
(u_j)^2\big)\int_{\Omega}\kern-3pt u_j\phi +\sum_{i\neq j}
2tG'\big(\int_{\Omega}\kern-3pt
(u_i)^2\big)\int_{\Omega}\kern-3pt u_i\phi+o(t) \\ &=&
t\int_{\Omega}\big(\nabla \widehat u_j\nabla\phi-
(\lambda_1(\omega_j)u_j-\sum_{i\neq j}
\lambda_1(\omega_i)u_i)\phi\big)+o(t)\;.
\end{array}
\endeq
For $t$ small enough we find $J(V)<J(U)$, a contradiction.

{\it The case $q>1$}. The idea of the proof is analogous, but a
new difficulty arises due to the fact that we can use only test
functions belonging to the Nehari manifold $\N_0$. Let us show the
new argument in the proof of the inequality 2. Assume by
contradiction that the assertion does not hold for a certain index
$i$ and  thus the existence of $0<\phi\in C^\infty_c$ such that
$$
\int_\Omega\Big[\nabla \widehat u_j\nabla\phi
-\Big(\lambda_1(\omega_j)u_j-\sum_{i\neq j}
\lambda_1(\omega_i)u_i\Big)\phi\Big]\,dx<0.
$$
We will obtain a contradiction constructing a $k$--uple in $\N_0$
that decreases the value of $c_q$. Let $\Lambda_j\widehat
u_j:=\lambda_j u_j-\sum_{i\neq j}\lambda_iu_i$ with $|\l_i-1|\leq
\delta$ for all $i$: if $\delta$ is small enough we can also
assume by continuity that
\neweq{aaus}
\oint\Big[\nabla \Lambda_j\widehat
u_j\nabla\phi-\Big(\lambda_1(\omega_j)\lambda_ju_j-\sum_{i\neq j}
\lambda_1(\omega_i)\lambda_iu_i\Big)\phi\Big]\;dx< 0.\endeq

By the inf--sup characterization of $c_q$ and by the behavior of
the function $J^*(\l u)$ for fixed  $u>0$, we can take $\delta$ so
small that
\neweq{aaus2}\nabla J^*((1-\delta)u_j)u_j>0,\hspace{1cm} \nabla
J^*((1+\delta)u_j)u_j<0, \hspace{1cm}\forall j.
\endeq
Let us fix $\bar{t}>0$ small and let us consider a $C^1$ function
$t:(\bbbr^+)^k\to \bbbr^+$ where
$t(\l_1,...,\l_k)=0$ if for at least one $j$ it happens
$|\l_j-1|\geq \delta$, and $t(\l_1,...,\l_k)=\bar{t}$ if
$|\l_j-1|\leq\delta/2$ for every $j$. Next we define
the
continuous map
$$\Phi(\l_1,...,\l_k)=\l_iu_i-\sum_{j\neq i}
\l_ju_j+t(\l_1,...,\l_k)\phi.$$ Note that $\Phi^-$ is a
positive function whose support is union of $k-1$ disjoint
connected components, each of them belonging to the support of
some $u_j$. Now we define the function
$\tilde{U}(\l_1,...,\l_k)=(\tilde{u}_1,...,\tilde{u}_k)$ as
$$
\tilde{u}_i=\left\{
\begin{array}{ll}
(\Lambda_j\widehat u_j+t(\l_1,...,\l_k)\phi)^+, &\mbox{if $i=j$},\\
(\Lambda_j\widehat u_j+t(\l_1,...,\l_k)\phi)^- \chi_{\{u_i>0\}}, &\mbox{if $i\neq j$}.
\end{array}\right.$$
 Let us compute $J(\tilde{U})$: in complete analogy with the calculations in \eq{co7} we have
$$
\begin{array}{c}
J(\tilde{U})=\int_\O\big( \frac{1}{2}|\nabla (\Lambda_j\widehat u_j)|^2+
\frac{t^2}{2}|\nabla \phi|^2\big)dx+
t\int_\O \nabla (\Lambda_j\widehat u_j)\nabla \phi\\
- G\big(\int_{\Omega}\kern-3pt((\lambda_j u_j+t\phi)^+)^2\big)+
\sum_{i\neq j}
G\big(\int_{\{u_i>0\}}\kern-3pt ((\Lambda_i\widehat
u_i+t\phi)^-)^2\big)\leq\\
\leq J(\l_1u_1,...,\l_k u_k)+\bar{t}\int_\O\big(\nabla(\Lambda_j\widehat u_j)\nabla
\phi-(\lambda_1(\omega_j)\lambda_ju_j-\sum_{i\neq j}
\lambda_1(\omega_i)\lambda_iu_i)\phi\big)
+o(\bar{t}).
\end{array}
$$ By \eq{aaus} and taking $\bar{t}$ small enough, this implies
$$J(\tilde{U}(\l_1,...,\l_k))<J(\l_1u_1,...,\l_ku_k)$$ if
$|\l_j-1|\leq \delta/2$ for every $j$.\\ Now, if $\bar{t}$ is
small, we can assume that \eq{aaus2} holds for
$\tilde{U}(\l_1,...,\l_k)$ instead of $(\l_1u_1,...,\l_k u_k)$.
Thus by continuity there exists $(\mu_1,...,\mu_k)$ such that
$|\mu_i-1|\leq \delta/2$ and
$$\nabla
J(\tilde{U}(\mu_1,...,\mu_k))\cdot\tilde{U}(\mu_1,...,\mu_k)=0$$
that means $\tilde{U}(\mu_1,...,\mu_k)\in{\mathcal N}_0$. But this is
in contradiction with the definition of $U$ as in Theorem
\ref{optimal} and the fact that
$$J(\tilde{U}(\mu_1,...,\mu_k))<J(\mu_1u_1,...,\mu_ku_k)\leq J(U)=
\inf_{V\in{\mathcal N}_0}J(V).$$
With this the proof of the inequality 2 is done.
Let us now briefly sketch the proof of the last inequality, namely
$$-\Delta u_i\leq \lambda_1(\omega_i)u_i$$
for all $i$. As usual
assume by contradiction the existence
of $\phi>0$, $\phi\in C^\infty_c(\Omega)$ such that
\neweq{con}
\oint\big[\nabla
(\l_iu_i)\nabla\phi-\lambda_1(\omega_i)\l_iu_i\phi\big]\;dx>0
\endeq
for all $\l_i$ such that $|\l_i-1|\leq \delta$, $\delta$ small
enough. As in the proof of the inequality 2, we can assume
$\delta$ small enough to satisfy \eq{aaus2}, and we consider the
function $t(\l_i)$ analogous to the one introduced therein. Then
we let $\Phi(\l_i):=\l_iu_i-t(\l_i)\phi$ and we define $\tilde{U}$
with components $\tilde{u}_i=\Phi^+$, $\tilde{u}_j=u_j$ if
$j=1,...,k-1$, $j\neq i$ and finally $\tilde{u}_k(x)=u_k(x)$ if
$x\in\{u_k>0\})$, $\tilde{u}_k(x)=\Phi^-$ if
$x\in\{u_i>0\}\cap
\{\Phi(x)<0\}$. By computing
$J(\tilde{U}(\l_i))$, taking into account \eq{con} and choosing
$\bar{t}$ small enough
 we obtain
$$
J(\tilde{U}(\l_i))<J(u_1,..,\l_iu_i,..u_k)
\hspace{1cm}|\l_i-1|\leq
\delta/2.$$ Now a contradiction with the properties of $U$ as in
Theorem \ref{optimal} can be obtained by arguing as in the final
step of the proof of the inequality 2.
\endproof

\subsection{Regularity results}\label{ultimasubsection}

Aim of this section is to prove Theorem \ref{extr_intr} in the
case $p\neq1$ (as already said, the case $p=1$ is treated in the
following section).

The differential inequalities obtained in Lemma \ref{super} allows
the application of the regularity theory developed by the authors
in \cite{ctv3}. To be more precise, consider a $k$--uple of
$H^1_0$ functions $(v_1,\dots,v_k)$; we set $\o_i:=\{v_i>0\}$,
$f_i(s):=\lambda_1(\omega_i)s$, $\widehat v_i:=v_i-\sum_{j\neq
i}v_j$ and
$$
\widehat f(x,\widehat v_i):=\lambda_1(\omega_i)v_i-\sum_{j\neq
i}\lambda_1(\omega_j)v_j.
$$
With these notations, if $(u_1,\dots,u_k)$ is as in Theorem
\ref{optimal}, then Lemma \ref{super} says that
$(u_1,\dots,u_k)\in{\mathcal S}$, where
$$
{\mathcal S}:=\left\{(v_1,\cdots,v_k)\in
(H^1(\Omega))^k:\,\begin{array}{l}
 v_i\geq0,\, v_i\cdot v_j=0\mbox{ if }i\neq j\\
 -\Delta v_i\leq f_i(x,v_i),\,-\Delta \widehat v_i\geq\widehat f(x,\widehat v_i)
\end{array}\right\}.
$$
This class of functions have been introduced by the authors in
\cite{ctv3}, where a number of qualitative properties for its
elements are obtained. We collect those properties
in the following theorem, referring to \cite{ctv3} for its proof.
We first need a definition:
\begin{defin}\label{def1}\rm
The multiplicity of a point $x\in\Omega$ with respect to the
$k$--uple $(v_1,\dots,v_k)$ is
$$
m(x)=\sharp\left\{i:meas\left(\left\{v_i>0\right\}\cap B_r(x)
\right)>0,\;\;\forall \,r>0\right\}\;.
$$
\end{defin}
\begin{teo}\label{main}
Let $(v_1,\dots,v_k)\in{\mathcal S}$, and let
$\omega_i:=\{v_i>0\}$. Then, $V:=\sum_{i=1}^k v_i$ verifies the
following properties.
\begin{enumerate}
\item The function $V$ is Lipschitz continuous in
the interior of $\O$; if $\partial\O$ is regular, then $V$ is
Lipschitz up to the boundary. In particular, $\omega_i=\{v_i>0\}$
is an open set.
\item Let $x\in\Omega$ such that $m(x)=2$. Then,
$$\lim_{y\to x\atop v_i(y)>0}\nabla v_i(y)=-\lim_{y\to x\atop v_j(y)>0}
\nabla v_j(y)\;.$$
\item In dimension $N=2$, the set $\{m(x)\geq3\}$ consists in a
finite numbers of points where $\nabla V$ is identically zero.
\item In dimension $N=2$, let $x\in\Omega$ such that $m(x)=h$. Then, $V$ admits a
local expansion around $x$ of the following form:
$$V(r,\theta)=r^{\frac{h}{2}}|\cos(\frac{h}{2}(\theta+\theta_0))|+o(r^{\frac{h}{2}})$$
as $r\to 0$, where $(r,\theta)$ denotes a system of polar
coordinates around $x$.
\item In dimension $N=2$, the set $\{m(x)=2\}$ consists in a
finite number of $C^1$--arcs ending either at points with higher
multiplicity, or at the boundary $\partial\O$.
\end{enumerate}
\end{teo}

By the above discussion, $(u_1,\dots,u_k)$ shares all these
properties. In particular, this implies that the partition
consisting of its supports, besides belonging to ${\mathcal
P}^*_k$, is an element of ${\mathcal P}_k$.

\begin{lemma}\label{connessione}
Let $(\omega_i)$ be the $k$--uple provided in Remark
\ref{deboleesistenza}. Then $\omega_i$ is open for every $i$.
Moreover we can assume, without loss of generality, that each
$\omega_i$ is connected (that is, $(\omega_i)\in {\mathcal P_k}$
and it is a solution of \eq{problema}).
\end{lemma}
\proof by Theorem \ref{optimal}, $\o_i=\{u_i>0\}$. Hence the
application of Theorem \ref{main},(1) provides that each $\o_i$ is
open. Assume that, for some $i$, $\o_i$ is not connected, and let
$\{\alpha_j\}_{j\in J}$ denote its connected (open) components. We
observe that, for every $M>0$, $\sharp\{j:\l_1(\a_j)\leq
M\}<\infty$ (indeed $\O$ is bounded and
$\lim_{|\o|\to0}\l_1(\o)=\infty$). Then
$\lambda_1(\omega_i)=\min_j\{\lambda_1(\alpha_j)\}=\lambda_1(\alpha_{h})$.
Replacing $u_i$ with $u_i|_{\alpha_{h}}$, we obtain a $k$--uple of
functions with open, connected supports, that again achieves
\eq{probweak}.
\endproof

{\bf Proof of Theorem \ref{extr_intr} (case $p\neq1$):} let
$(u_1,\dots,u_k)$ the functions provided by Theorem \ref{optimal}
and let $\omega_i=\{u_i>0\}$. Clearly, $u_i=a_i\phi_i$, where
$\phi_i$ denotes the positive eigenfunction associated to
$\lambda_1(\omega_i)$, normalized in $L^2$, and
$a_i=\lambda_1^{1/2(q-1)}(\omega_i)$. As shown in Lemma
\ref{super}, each $u_i$ satisfies the required differential
inequalities. Therefore, Theorem \ref{main} applies, implying that
each $\o_i$ is open. Finally, also Lemma \ref{connessione}
applies, and the theorem follows.\endproof

A remarkable consequence of the above results is the equivalence between the original problem \eq{problema} and the relaxed one \eq{problema*}:
\neweq{forte=debole}
\inf_{(\omega_i)\in{\mathcal P}_k}\sum_{i=1}^k
(\lambda_1(\omega_i))^p\equiv\inf_{(\omega_i)\in{\mathcal
P}^*_k}\sum_{i=1}^k (\lambda_1(\omega_i))^p.
\endeq

Up to now, this is true only when $p\neq1$. The
discussion in the following section will trivially imply that it
holds also in the case $p=1$.

\begin{rem}\label{gener}{\rm
Let us conclude this section with a remark about the generality of
the theory so far developed. Actually, the procedure leading to
the proof of Theorem \ref{optimal} and, consequently, to Theorem
\ref{extr_intr}, can be trivially adapted to study, for instance,
nonhomogeneous optimal partition problems. Namely, let $m_i,n_i\in
L^\infty$ such that $\inf_{x\in\Omega}\{m_i(x),n_i(x)\}>0$;
finally let $a_i\in \bbbr^+$ and  $q\in\bbbr$, $r\in\bbbr$ such
that $q>1$ and $r\geq 2$. By defining the first
weighted--eigenvalue as
$$\lambda_1(\omega_i)=\min_{u\in H^1_0(\omega_i)\atop u\not\equiv 0}
\frac{\int_{\omega_i}n_i(x) |\nabla
u(x)|^r dx}{\int_{\omega_i}m_i(x)|u(x)|^r dx},$$ we consider the problem of
finding a partition of $\O$ in $k$ open sets that achieves
\neweq{co1}c_q\equiv\inf_{(\omega_i)\in{\mathcal P}_k}
\frac{q-1}{rq}\sum_{i=1}^k
[a_i\lambda_1(\omega_i)]^\frac{q}{q-1}.\endeq Then, with
obvious changes in the functional setting, namely by redefining
$$
\dis J_i(u)=\dis\frac{a_i}{r}\int_\O n_i(x)|\nabla u(x)|^r dx -
\dis\frac{1}{rq}\left(\int_\O m_i(x)|u(x)|^r dx\right)^q
$$
the whole procedure applies to problem
\eq{co1}.
Note that this includes the remarkable case of the $r$--laplacian, and will be crucial
in connection with the analysis of the Fu\v{c}\'\i k spectrum developed in the last section.}
\end{rem}

\section{The limiting cases $p=1$, $p=\infty$}\label{p=1}

In this section we study the asymptotic behavior of the solutions
to the  problem \eq{ck} both as $p\to 1$ and $p\to \infty$. This
analysis will provide existence and regularity results, analogous
to those obtained in the previous section, for two remarkable
optimal partition problems that we cannot directly treat with the
above techniques. Let us start our description by observing that
the eigenvalues corresponding to \eqref{problema} are uniformly
bounded in $p$:
\begin{lemma}\label{mM}
Let $p>0$, $p\neq 1$ and let $(\omega_{1,p},...\omega_{k,p})\in
{\mathcal P}_k$ achieving \eqref{problema}. Then there exist
$0<m<M<\infty$ such that
$$m\leq \lambda_1(\omega_{i,p})\leq M\hspace{1cm}\forall
i,\;\forall p.$$
\end{lemma}
\proof the bound from below depends on the monotonicity of the
first eigenvalue with respect to the inclusion
$$
\omega\subset\Omega \Longrightarrow\lambda_1(\omega)\geq \lambda_1(\Omega).
$$
Hence it suffices letting $m=\lambda_1(\Omega)$.
The bound from above simply follows by the minimality of the optimal
partition.
\endproof
In the following we shall denote $\lambda_{i,p}:=\lambda_1(\omega_{i,p}).$
Let us also recall that, if $u_{i,p}$ are the eigenfunctions
associated to $c_q$, then by Theorem \ref{optimal}
\neweq{limdeb}
\left(\int_\O u_{i,p}^2\right)^{q-1}=\lambda_{i,p}=
\left(\int_\O |\nabla u_{i,p}|^2\right)^\frac{q-1}{q},
\endeq
for all $i$.

\subsection{Proof of Theorem \ref{extr_intr} (case $p\neq1$):}
let
\neweq{problema_L_1} c_1= \inf_{(\omega_i)\in{\mathcal
P}_k}\sum_{i=1}^k \lambda_1(\omega_i),
 \endeq
and let $p\to1$, hence $q\to\infty$: since $\lambda_{i,p}$ is
uniformly bounded in $p$, we obtain by \eq{limdeb} that $u_{i,p}$
is bounded in $H^1$. Therefore, there exist $u_i\in H^1$ such that
(up to a subsequence) $u_{i,p}\rightharpoonup u_i$ weakly in
$H^1$. Since the convergence is also almost everywhere, then,
calling $\omega_i:=\{u_i>0\}$, we have that the $(\omega_i)$'s are
disjoint. We claim that $(\omega_i)$ is a solution of
\eq{problema_L_1} and that the corresponding $u_i$ satisfy
suitable extremality conditions.

To start with, let us observe that, by virtue of Lemma \ref{mM},
there exist $\mu_i\in[m,M]$ such that, up to a subsequence,
$\lim_{p\to 1}\lambda_{i,p}=\mu_i.$ By weak convergence, the
differential inequalities for $u_{i,p}$ pass to the limit, namely
$$-\Delta u_i\leq \mu_i u_i \hspace{1cm}\mbox{ and }\hspace{1cm}
-\Delta\widehat u_i\geq \mu_iu_i-\sum_{j\neq i} \mu_ju_j.$$
This allows to prove that the weak convergence is indeed strong.
To this aim, let us test the inequality $-\Delta\widehat u_i\geq
\mu_iu_i-\sum_{j\neq i} \mu_ju_j$ with $u_i$: then
$$
\int_\O |\nabla u_{i}|^2\geq \mu_i\int_\O u_i^2
$$
On the other side, by $-\Delta u_{i,q}\leq \lambda_{i,q} u_{i,q}$
tested with $u_{i,q}$ it holds
$$
\int_\O |\nabla u_{i,q}|^2\leq \lambda_{i,q}\int_\O u^2_{i,q}.
$$
By gluing the two previous inequality when passing to the limit we
obtain
$$\int_\O |\nabla u_{i}|^2\geq \limsup\int_\O |\nabla u_{i,q}|^2.$$
This finally provides $u_{i,n}\to u_{i}$ in $H^1$. Furthermore, by the
variational characterization of the first eigenvalue we have
$\mu_i\equiv \lambda_1(\omega_i)$.

As a consequence of this analysis we have that $c_q\to c_1$ and
that $(\omega_i)$ achieves $c_1$; furthermore it holds $$-\Delta
u_i\leq \lambda_1(\omega_i) u_i\hspace{1cm}\mbox{ and
}\hspace{1cm} -\Delta\widehat u_i\geq
\lambda_1(\omega_i)u_i-\sum_{j\neq i}\lambda_1(\omega_j)u_j,$$ as
required. But now we are in a position to apply the already
mentioned regularity theory (Theorem \ref{main}, Lemma
\ref{connessione}), providing $(\omega_i)\in {\mathcal P}_k$ and
concluding the proof.
\endproof

\subsection{The case $p=\infty$.}\label{secondo.autovalore}

Let $p\to\infty$ (hence  $q\to 1$). By virtue of the basic
property
$$\lim_{p\to\infty}\left(\sum_{i=1}^k
a_i^p\right)^\frac{1}{p}=\max\{a_1,...,a_k\},$$ where $a_i$ are
positive numbers, we shall succeed in recovering our existence and regularity results for a partition achieving
\neweq{problema_L_2}
c_\infty:=\inf_{(\omega_i)\in{\mathcal P}_k}\max_{i=1,...,k}
\lambda_1(\omega_i).
\endeq
Indeed we are going to prove
\begin{teo}
There exists $U\in{\mathcal S}$ such that $(\{u_1>0\},...,\{u_k>0\})$ achieves the value
\eq{problema_L_2}.
\end{teo}
\proof let $p\to\infty$; note that in this case we do not know if
$u_{i,q}$ is $H^1$--bounded. But if we define
$$v_{i,q}=\frac{u_{i,q}}{\left(\sum\lambda_{i,q}\right)^\frac{q}{q-1}},$$
then \eq{limdeb} ensures $\|v_{i,q}\|_{H^1}\leq 1$ for all $i$ and
$q$. Hence $v_{i,q}$ admit a weak limit (up to a subsequence) and
all the analysis developed in the previous case still holds. In
particular, if we call $u_i$ the $H^1$ limit of $v_{i,q}$, again
the extremality conditions hold true and consequently
$(\{u_i>0\})$ belongs to ${\mathcal P}_k$. Let us now prove that
$(\{u_i>0\})$ is indeed a solution for \eq{problema_L_2}. By the
strong convergence of the $v_{i,q}$'s and the above mentioned
basic property we know that
$$c_q^\frac{1}{p}=\left(\sum_{i=1}^k
\lambda_{i,q}^p\right)^\frac{1}{p}\to
\max\{\lambda_1(\{u_1>0\}),...,\lambda_1(\{u_k>0\})\}:=M_k,$$
hence it is enough to prove that $M_k\equiv c_\infty$. To this aim
let $(\omega_1,...,\omega_k)\in {\mathcal P}_k$ a $k$--uple of
disjoint sets achieving $c_\infty$. Then by definition it holds
$$
c_q^\frac{1}{p}\leq
\left(\sum_{i=1}\lambda_1(\omega_i)^p\right)^\frac{1}{p}\to
\max\{\lambda_1(\omega_1),...,\lambda_1(\omega_k)\}\equiv
c_\infty.
$$
This implies
$$
M_k=\lim_{p\to\infty}c_q^\frac{1}{p}\leq c_\infty.
$$
By the minimality of $c_\infty$ the opposite inequality
$c_\infty\leq M_k$ is immediate, thus $M_k\equiv c_\infty$ as
claimed.
\endproof
\section{The first curve of the Fu\v{c}\'\i k spectrum}
Let us consider the problem
\neweq{fuc}-\Delta u=\lambda u^+-\mu u^-\endeq
in $\O\subset \bbbr^N$ with boundary condition $u=0$ on
$\partial\O$; the Fu\v{c}\'\i k spectrum of $-\Delta $ on $H^1_0(\O)$ is defined as
$${\mathcal F}=\{(\lambda, \mu)\in\bbbr^2\;:\;\mbox{ problem  \eq{fuc}
has a non--trivial (weak) solution }\}\;.$$ As already discussed,
this object has been argument of a quite extensive literature
devoted to study its structure and its connections with the
solvability of nonlinear related problems (see references in the
introduction). In particular, it is known by \cite{dfg} that,
besides the pairs of equal eigenvalues and the semi--lines
$(\lambda_1(\O), t)$, $(t, \lambda_1(\O))$ for all $t>0$,
${\mathcal F}$ contains a first nontrivial curve $C_1$ through
$(\lambda_2(\O), \lambda_2(\O))$, which extends to infinity. The
objective of this section is to give a new description of $C_1$ in
the variational setting developed in the present paper. To this
aim, let $r$ be a positive number; we introduce the problem of
finding a partition achieving
\neweq{c.r}
 c(r):=\inf_{(\omega_i)\in{\mathcal P}_2}
 \max\{r\lambda_1(\omega_1),\lambda_1(\omega_2)\}.
\endeq
Our characterization is given by Theorem \ref{fuc_intr} as stated
in the introduction. We recall the result:

{\bf Theorem \ref{fuc_intr}.}$\;$
{\it For all $r>0$, there exists $u\in H^1_0(\Omega)$ such that
$(\{u^+>0\},\{u^->0\})$ achieves $c(r)$. Furthermore, the pair
\neweq{lm}(\lambda_1(\{u^+>0\}),\lambda_1(\{u^->0\}))=(r^{-1}c(r),
c(r))\endeq belongs to the Fu\v{c}\'\i k spectrum ${\mathcal F}$
 and it represent the first (nontrivial) intersection between ${\mathcal F}$
and the line of slope $r$.}

As a consequence we have a
variational characterization of the second eigenvalue of the Laplacian in $H^1_0(\Omega)$:
\begin{coro}
$$
 \lambda_2(\Omega)\equiv\inf_{\omega\in\Omega}
 \max\{\lambda_1(\omega),\lambda_1(\Omega\setminus\omega)\}.
$$
\end{coro}
\proof note that, for the choice $r=1$, in view of \eq{remark}, it
holds that $c(1)$ is an eigenvalue corresponding to a
sign--changing eigenfunction. Hence $c(1)\geq\lambda_2(\O)$. On
the other hand, it follows from the property of $c(1)$ of being
the {\it first}  intersection with the  Fu\v{c}\'\i k spectrum
${\mathcal F}$ and the well--known fact
$(\lambda_2(\O),\lambda_2(\O))\in{\mathcal F}$ that
$c(1)=\lambda_2(\O)$. From this the thesis follows.\endproof

{\bf Proof of Theorem \ref{fuc_intr}:}  let $q>1$ be fixed: by
Remark \ref{gener}, with the choice $k=2$, $m_i(x)=n_i(x)=1$ and
$a_1=r$, $a_2=1$,  we immediately obtain the existence of a pair
$(u_{1,q},u_{2,q})$ whose supports $(\omega_{1,q},\omega_{2,q})$
achieve the value
$$
c_q\equiv\inf_{(\omega_i)\in{\mathcal P_2}} \frac{q-1}{2q}
\Big((r\lambda_1(\omega_{1}))^\frac{q}{q-1}+\lambda_1(\omega_2)^\frac{q}{q-1}\Big).
$$
Now we choose a sequence $q\to 1$ and we follow the limiting
procedure in Section \ref{secondo.autovalore} and the arguments
therein. We thus obtain, when passing to the limit as $q\to 1$,
the existence of a pair of $H^1_0(\O)$--functions $(u_1,u_2)$ with
the following properties. First, $U=(u_1,u_2)\in {\mathcal S}$,
and $(\{u_1>0\},\{u_2>0\})$ achieves the value
$$\inf_{(\omega_i)\in{\mathcal
P_2}}\max\{r\lambda_1(\omega_1),\lambda_1(\omega_2)\}.
$$
Let us now  define the function $u:=u_1-u_2$, that is, $u_1=u^+$
and $u_2=u^-$. Since  $U\in{\mathcal S}$, then
properties  (1) and (2) of Theorem \ref{main} hold. Hence
$u$ is regular and it is a nontrivial solution of
$$
-\Delta u= \lambda_1(\{u^+>0\})u^+-\lambda_1(\{u^->0\})u^-
$$
on $\Omega$. We have
\neweq{remark}
r\lambda_1(\{u^+>0\})=\lambda_1(\{u^->0\})=c(r).
\endeq
Indeed, assume by contradiction that $c(r)=r\lambda_1(\{u^+>0\})$
and $r\lambda_1(\{u^+>0\})-\lambda_1(\{u^->0\})=m>0$. Let $x\in
\partial\{u^+>0\}\cap\partial\{u^->0\}$. Since $u$ is regular, we
can choose $\rho>0$ small enough in such a way that, by the
monotonicity and continuity properties of the first eigenvalue,
there holds
$$
r\lambda_1(\{u^+>0\})>r\lambda_1(\{u^+>0\}\cup B(x,\rho))>
r\lambda_1(\{u^+>0\})-\frac{m}{4}
$$
and
$$
\lambda_1(\{u^->0\})<\lambda_1(\{u^->0\}\setminus B(x,\rho))<
\lambda_1(\{u^->0\})+\frac{m}{4}.
$$
In this way we have a new partition $(\{u^+>0\}\cup B(x,\rho),
\{u^->0\}\setminus B(x,\rho))$ which lowers the value $c(r)$, a
contradiction. As a consequence, we finally obtain that the pair
$(r^{-1}c(r),c(r))$ belongs to $\mathcal F$.

We are left to prove that \eq{lm} is in fact the {\it first}
nontrivial intersection of the spectrum with the line of slope
$r$. First we observe that,  by the monotonicity of the first
eigenvalue with respect to the inclusion, $c(r)>\lambda_1(\Omega)$
for every $r>0$. Assume by contradiction the existence of a pair
$(r^{-1}\mu,\mu)\in {\mathcal F}$ such that $\mu<c(r)$. This means
the existence of $v\in H^1_0$ with $v^\pm\not\equiv 0$ such that
$-\Delta v= r^{-1}\mu v^+-\mu v^-$: testing the equation with
$v^\pm$ we have
$$r^{-1}\mu=r^{-1}\frac{\int_\O |\nabla v^-|^2}{\int_\O |v^-|^2}=
\frac{\int_\O |\nabla v^+|^2}{\int_\O |v^+|^2}.$$
Since $(\{v^+>0\},\{v^->0\})$ is an admissible partition of $\O$, it must hold
$$c(r)\leq \max \left\{ r\frac{\int_\O |\nabla v^+|^2}{\int_\O |v^+|^2},
\frac{\int_\O |\nabla v^-|^2}{\int_\O |v^-|^2} \right\}=\mu,$$ a
contradiction.
\endproof
Hence, by denoting
$$C_1:=\{(r^{-1}c(r),c(r))\in\bbbr^2,\;r>0\},$$
we have that $C_1$ is indeed the first nontrivial curve of the Fu\v{c}\'\i k spectrum.
We wish to emphasize that our variational characterization of
$C_1$  immediately provides the main feature of the first curve
and of the eigenfunctions associated to each element of $C_1$ (see
\cite{dfg}). In particular, just by reading the definition of $c(r)$
we can prove the following
\begin{prop}
$\;$
\begin{enumerate}
\item[(a)] $C_1$ is a continuous and strictly decreasing curve,
symmetric with respect to the diagonal.
\item[(b)] $C_1\subset \{(x,y)\in\bbbr^2\; :\; x>\lambda_1(\O), y>\lambda_1(\O)\}$
and it is asymptotic to the lines $\lambda_1(\Omega)\times \bbbr$
and $\bbbr\times\lambda_1(\Omega)$.
\item[(c)] Any eigenfunction associated to $(x,y)\in C_1$ admits exactly two nodal domains
(Courant nodal domain theorem).
\end{enumerate}
\end{prop}
\proof (a) The symmetry of $C_1$  can be derived  by
interchanging the role of $\omega_1$ and $\omega_2$. The
continuity it's immediate by the definition of $c(r)$ as in
\eq{c.r}; the monotonicity of the curve is equivalent to the fact
that $r_1>r_2\Rightarrow c(r_1)>c(r_2)$: but this directly follows
once again by the definition of $c(r)$.\\
(b) The first part is given by the fact that $c(r)>\lambda_1(\O)$.
Proving  the existence of the asymptotes is equivalent to show
that $c(r)\to\lambda_1(\Omega)$ as $r\to 0$. To this aim, let
$\eps>0$ be fixed and consider a small ball of radius $\rho$, such
that $\lambda_1( \O\setminus B(x,\rho))\leq\lambda_1(\O)+\eps$.
Let us consider the partition made up by  $\omega_1=B(x,\rho)$ and
$\omega_2=\O\setminus B(x,\rho)$, then choose $r>0$ small enough
in such a way that $r\lambda_1(B(x,\rho))< \lambda_1(\O)+\eps$: it turns out
that
$c(r)\leq\max\{r\lambda_1(\omega_1),\lambda_1(\omega_2)\}\leq\lambda_1(\O)+\eps$.
We have thus proved that $\forall \eps>0$ there exists $r=r_\eps>0$
such that $\lambda_1(\O)<c(r_\eps)<\lambda_1(\O)+\eps$, concluding
the proof.\\
 (c) The nodal property its already true by definition
of $c(r)$ and the procedure of partitioning $\Omega$ exactly in
two connected subsets (see Lemma \ref{connessione}).
\endproof

\subsection{Further results.} In this section let us develop some extensions
of the previous techniques. The first applies to the search of further elements of the
Fu\v{c}\'\i k spectrum. Then, we show how to recover the case of more general operators.

For $1\leq h\leq k$, let us define the numbers
 $$c_{h,k}(r):=\inf_{(\omega_i)\in{\mathcal P}_k}
 \max\{r\lambda_1(\omega_1),\dots,r\lambda_1(\omega_h),\lambda_1(\omega_{h+1}),\dots,\lambda_1(\omega_k)\}\;.
 $$
 We know by all the previous discussion that the infima above are attained and that a suitable
 choice of the eigenfunctions $u_i$ corresponding to the optimal partition $\omega_i$,
 satisfy the extremality conditions stated in
 Theorem \ref{extr_intr}.  Moreover we have that
 $$c_{h,k}(r)\equiv \lambda_1(\omega_j)=r\lambda_1(\omega_i),\qquad  1\leq i\leq h,\,
 j\geq h+1.$$
 Now assume that the interfaces $\Gamma_{i,j}=\partial \omega_i\cap \partial\omega_j$ consist only of points of multiplicity two and moreover that $\Gamma_{i,j}=\emptyset$ unless $i\in\{1,\dots,h\}$ and
 $j\in\{h+1,\dots,k\}$. Let us consider the function
 $$u=u_1-u_{h+1}+u_2-u_{h+2}+\dots.$$
 Then, by property (2) in Theorem \ref{main}, we have that $u$ is regular
 and thus it is a solution of the equation
 $$-\Delta u=r^{-1}c_{h,k}(r)u^+-c_{h,k}(r)u^-,$$
 in the whole of $\O$.
 Therefore, under suitable topological assumptions,  we have proved that
   $$(r^{-1}c_{h,k}(r),c_{h,k}(r))\in{\mathcal F},$$
  providing a new nontrivial element of the spectrum.
  This may happen in several practical situations, as shown by the numerical experiments
  in \cite{hr}. For instance, this is
 always the case when working in 1--dimensional domains. This actually proves
 the existence of
 a sequence of curves of the Fu\v{c}\'\i k spectrum, as stated in
 Theorem \ref{altre-curve}.

 \medskip

 Let us conclude by pointing out some possible generalizations
 of the above results. Actually, thanks to the discussion developed in Remark
 \ref{gener}, we already know that the abstract setting leading to the proof of Theorem \ref{main},
 applies to more general problems. As a consequence, we can prove results analogous to Theorem
 \ref{fuc_intr} and \ref{altre-curve}, which describe nontrivial elements of some possible
 generalization of the spectrum. In particular, we can characterize a first nontrivial curve
 of elements in the spectrum of  the $p$--Laplacian (see \cite{dg1} for the definition and
 a comparison) just replacing the notion of first eigenvalue with the one related to the
  new operator, namely
  $$\lambda_1(\omega):=
\min_{u\in W^{1,p}_0(\omega)\atop u\not\equiv 0} \frac{\int_\omega |\nabla
u(x)|^pdx}{\int_\omega|u(x)|^pdx}.$$
 Another interesting application consists in the characterization of the first curve
 of elements for a generalized notion of
 spectrum in presence of positive weights $p,q$. Namely, Theorem \ref{fuc_intr}
applies to describe the set of $(\lambda,\mu)$ such that
 $$-\Delta u=\lambda p(x)u^+-\mu q(x)u^-$$
 has a nontrivial solution.
 In this case we have the natural replace of the definition of $\lambda_1$ with the
 corresponding weighted ones
$$\lambda_1(\omega_1)=\min_{u\in H^1_0(\omega_1)\atop u\not\equiv 0} \frac{\int_{\omega_1} |\nabla
u(x)|^2
dx}{\int_{\omega_1}p(x)|u(x)|^2dx},\quad\qquad\lambda_1(\omega_2)=
\min_{u\in H^1_0(\omega_2)\atop u\not\equiv0}
\frac{\int_{\omega_2} |\nabla u(x)|^2
dx}{\int_{\omega_2}q(x)|u(x)|^2dx}.$$

 Finally, let us only mention that the whole theory
 can be easily modified in order to apply to general boundary conditions
besides the Dirichlet case.


\section{Monotonicity Formulae}

Consider the general problem of minimizing
\neweq{valore}
\beta(k,N):= \inf_{{\mathcal P}(k,N)}\frac{2}{k}\sum_{i=1}^k , \sqrt{\lambda_1(\omega_i)}
\endeq
where $S^{N-1}$ denotes the boundary of the unit ball in $\bbbr^N$
and
$$
{\mathcal P}(k,N):=\left\{(\omega_1,\dots,\omega_k)\subset
S^{N-1}:\;\omega_i\mbox{ is open and
connected},\;\omega_i\cap\omega_j=\emptyset\;{\rm if}\;i\neq
j\right\}.
$$
In this section we are concerned with the properties of
$\beta(k,N)$ and with its relation with some extensions of the
monotonicity formula. First of all, it can be proved that
$\beta(k,N)$ is achieved by a partition containing only open and
connected sets of $S^{N-1}$. This directly comes by Remark
\ref{deboleesistenza}, Theorem \ref{main}, and Lemma
\ref{connessione}, where the results are obtained for partitions
of domains in $\bbbr^N$. The proof of this fact for partitions of
$S^{N-1}$ can be recovered (possibly through local charts) in a
straightforward way.

Let us now concentrate on the value of $\beta$: when there are
only two parts, the optimal partition is achieved by the
equator--cut sphere (see \cite{sp}) and hence
$$
\beta(2,N)=N,
$$
(thus, in particular, for $k=2$ and $N=2$ our Lemma
\ref{monot.generale_intr} exactly gives the result in \cite{cks}).
When $k\geq3$ the only exact value of $\beta$ we can give refers
to the dimension $N=2$ and reads
$$
\beta(k,2)=k,
$$
as can be found in \cite{ctv2}. Nevertheless, for $k\geq 3$, we
are going prove that
$$
\beta(k,N)>N
$$
in any dimension larger than 2, as a consequence of the
monotonicity of $\beta$ as a function of $k$:
\begin{prop}\label{beta.cresce} The function
$\beta(\cdot,N):\bbbn\to\bbbr^+$ is non decreasing. Moreover,
$\beta(k,N)>\beta(2,N)$ for $k\geq3$.
\end{prop}
\proof let $N\geq2$ be fixed and let
$(\Omega_1,\dots,\Omega_{k+1})\in{\mathcal P}(k+1,N)$ be a
partition of $S^{N-1}$ which achieves $\beta(k+1,N)$. Let us
assume, to fix the ideas, that $\l_1(\Omega_{k+1})\geq\l_1
(\Omega_i)$, $i=1,...,k$. If we consider $(\O_1,\dots,\O_k)$ as an
element of ${\mathcal P}(k,N)$ we easily obtain
$$
\beta(k+1,N)=\frac{1}{k+1}\sum_{i=1}^{k+1}
\sqrt{\lambda_1(\Omega_i)}\geq\frac{1}{k}\sum_{i=1}^{k}
\sqrt{\lambda_1(\Omega_i)}\geq\beta(k,N),
$$
and the equality holds iff $\l_1(\O_i)=\l_1(\O_{k+1})$ for every
$i$ and $(\O_1,\dots,\O_k)$ achieves $\beta(k,N)$. As a first
consequence, we obtain the weak monotonicity of $\beta(\cdot,N)$.

To conclude the proof of the lemma, we will show that $ \beta(3,N)
>\beta(2,N)$. Assume by contradiction that
$\beta(3,N)=\beta(2,N)$. By the above considerations (in the case
$k=2$), we obtain that $\l_1(\O_1)=\l_1(\O_2) =\l_1(\O_3)=:\l_2$
and that $(\O_1,\O_2)$ achieves $\beta(2,N)$. Let
$(u_1,u_2)\in{\mathcal S}$ denote the associated eigenfunctions.
Then, by definition of ${\mathcal S}$, we obtain both
$-\Delta\widehat u_1\geq\lambda_1(\Omega_1)
u_1-\lambda_1(\Omega_2)u_2$ and $-\Delta\widehat
u_2\geq\lambda_1(\Omega_2) u_2-\lambda_1(\Omega_1)u_1$, that is,
$-\Delta(u_1-u_2)=\lambda_2(u_1-u_2)$ on $S^{N-1}$; but
$u_1-u_2\equiv0$ on $\O_3$, in contradiction with the well known
properties of the eigenfunctions of the Laplace operator (unique
continuation property).
\endproof

As we said, the function $\beta$ naturally appears when trying to
extend a variant of the monotonicity formula to the case of many
subharmonic densities. In this perspective we are going to prove
Lemma \ref{monot.generale_intr}.

\medskip
{\bf Proof of Lemma \ref{monot.generale_intr}:} the idea of the
proof consists in showing that $\Phi'(r)\geq0$ for every $r$. Let us start with some estimates.
First,
since each $w_i$ is positive and $-\Delta w_i\leq0$,
testing with $w_i$ on the sphere $B(x_0,r)=:B_r$ (with
$r\leq\bar{r}$) we obtain, for every $i$:
$$
\int_{B_r}|\nabla w_i|^2\leq\int_{\partial B_r}w_i
\frac{\partial}{\partial n} w_i.
$$
Let $\nabla_T w_i:=\nabla w_i-n\partial_n w_i$ be the tangential
component of the gradient of $w_i$. We apply the H\"older inequality
to the previous equation, then we multiply and divide by the
$L^2$--norm of $\nabla_T w_i$, and finally we use the Young
inequality. We have
\neweq{prelimmonot}
\begin{array}{cl}
 \int_{B_r}|\nabla w_i|^2 &\leq  (\int_{\partial B_r}
 w_i^2)^{1/2}(\int_{\partial B_r}(\partial_n w_i)^2)^{1/2}\leq\\
 &\leq  (\int_{\partial B_r}|\nabla_T w_i|^2)^{1/2}(\int_{
 \partial B_r}(\partial_n w_i)^2)^{1/2}\cdot\frac{(\int_{\partial
 B_r}w_i^2)^{1/2}}{(\int_{\partial B_r}|\nabla_T w_i|^2)^{1/2}}\leq\\
 &\leq  \half\left[\int_{\partial B_r}|\nabla_T w_i|^2+\int_{
 \partial B_r}(\partial_n w_i)^2\right]\cdot\left(\frac{\int_{\partial
 B_r}w_i^2}{\int_{\partial B_r}|\nabla_T w_i|^2}\right)^{1/2}\leq\\
 &\leq  \half\int_{\partial B_r}|\nabla w_i|^2\cdot\left(\frac{\int_{\partial
 B_r}w_i^2}{\int_{\partial B_r}|\nabla_T w_i|^2}\right)^{1/2}.
\end{array}
\endeq
Now let $v_i^{(r)}:S^{N-1}\to\bbbr$ be defined as
$v_i^{(r)}(\xi):=w_i(x_0+r\xi)$, in such a way that $\nabla
v_i^{(r)}(\xi) =r^2\nabla_T w_i(x_0+r\xi)$. By the previous
inequality we obtain
$$
\frac{\int_{\partial B_r}|\nabla w_i|^2}{\int_{B_r} |\nabla
w_i|^2}\geq2\left(\frac{\int_{\partial B_r}|\nabla_T
w_i|^2}{\int_{\partial
B_r}w_i^2}\right)^{1/2}\geq2\left(\frac{r^{-2}\int_{S^{N-1}}|\nabla
v_i^{(r)}|^2}{\int_{S^{N-1}}(v_i^{(r)})^2}\right)^{1/2}\geq\frac2r\sqrt{
\l_1(\{v_i^{(r)}>0\})}.
$$
Since $w_iw_j=0$ a.e., the supports of the $v_i^{(r)}$'s
constitute a partition of $S^{N-1}$. Therefore, summing up on $i$
the previous inequality and recalling the definition of $\beta$
(and also \eq{forte=debole}), we finally have
\neweq{bordoBsuB1}
\sum_{i=1}^h \frac{\int_{\partial B_r}|\nabla w_i|^2}{\int_{B_r}
|\nabla w_i|^2}\geq\frac{h}{r}\beta(h,N).
\endeq
Now we are ready to prove the lemma: by computing $\Phi'(r)$ we
obtain
$$
\begin{array}{rcl}
 \dis\Phi'(r)&=&\dis-\frac{h\beta(h,N)}{r^{h\beta (h,N)+1}}\prod_{
 i=1}^{h}\int_{B_r}|\nabla w_i|^2+\frac{1}{r^{h\beta(h,N)}}\sum_{
 i=1}^h\left[\left(\prod_{j\neq i}\int_{B_r}|\nabla w_j|^2\right)
 \int_{\partial B_r}|\nabla w_i|^2\right]=\\
 &=&\dis\Phi(r)\left(-\frac{h\beta(h,N)}{r}+\sum_{i=1}^h
 \frac{\int_{\partial B_r}|\nabla w_i|^2}{\int_{B_r}|\nabla
 w_i|^2}\right).
\end{array}
$$
Replacing \eq{bordoBsuB1} in the previous equation, the lemma
follows.
\endproof

\begin{rem}
The argument above  shows that the function $\Phi$ is in fact strictly increasing, except in the case
when $w_i(r,\theta)=r^{\alpha}\phi_i$, where $\alpha=\beta(k,N)-N+1$ and the $\phi_i$'s are the first eigenfunctions of
the Laplace--Beltrami operator on the unit sphere, associated to the optimal partition \eq{valore}.
\end{rem}

With similar ideas we can prove also Lemma
\ref{monot.formula.N_intr}.

{\bf Proof of Lemma \ref{monot.formula.N_intr}:} we follow the
outline of the proof of Lemma \ref{monot.generale_intr}: we test
the equation with $u_i$ and, after some calculations, we obtain
the counterpart of equation \eq{prelimmonot}, that is
$$
\int_{B_r}\big[|\nabla u_i|^2+u_i^2\sum_{j\neq
i}a_{ij}u_j\big]\leq\half \int_{\partial B_r}\big[|\nabla
u_i|^2+u_i^2\sum_{j\neq
i}a_{ij}u_j\big]\cdot\left(\frac{\int_{\partial
B_r}u_i^2}{\int_{\partial B_r}\big[|\nabla u_i|^2+u_i^2\sum_{j\neq
i} a_{ij}u_j\big]}\right)^{1/2}.
$$
As in the proof of that lemma, we let $v_i^{(r)}:S^{N-1}\to\bbbr$
be defined as $v_i^{(r)}(\xi):=u_i(r\xi)$, and again $\nabla
v_i^{(r)}(\xi) =r^2\nabla_T u_i(r\xi)$, providing
$$
 \frac{\int_{\partial B_r}|\nabla u_i|^2+u_i^2\sum_{j\neq
 i}a_{ij}u_j}{\int_{B_r}|\nabla u_i|^2+u_i^2\sum_{j\neq
 i}a_{ij}u_j}\geq2\left(\frac{\int_{\partial B_r}|\nabla
 u_i|^2+u_i^2\sum_{j\neq i} a_{ij}u_j}{\int_{\partial
 B_r}u_i^2}\right)^{1/2}\geq
 \frac{2}{r}\sqrt{\Lambda_i(r)}
$$
where
$$
\Lambda_i(r)=\frac{\int_{S^{N-1}}|\nabla v_i^{(r)}|^2+r^2
(v_i^{(r)})^2\sum_{j\neq
i}a_{ij}v_j^{(r)}}{\int_{S^{N-1}}(v_i^{(r)})^2}.
$$
By computing $\Phi'$ as in the proof of Lemma
\ref{monot.generale_intr}, and taking into account the previous
calculations, we have
$$
\Phi'(r)\geq \Phi(r)\left(-\frac{hh'}{r}+\frac{2}{r} \sum_{i=1}^h
\sqrt{\Lambda_i(r)}\right).
$$
Hence we are lead to prove that there exists $r_0$ sufficiently
large, such that
\neweq{tesi}
\Phi'(r)\geq0\qquad \forall r\geq r_0.
\endeq
Observe that here we can not conclude as in the proof
of Lemma \ref{monot.generale_intr}: indeed, the supports of the
functions $v_i^{(r)}$'s are not mutually disjoint, and thus we can
not compare the value of $\sum \sqrt{\Lambda_i(r)}$ with
$\beta(h,N)$. To overcome this problem, we will let $r\to\infty$,
proving the convergence of (suitable multiples of) the
$v_i^{(r)}$'s to a $k$--uple of functions on $S^{N-1}$ having
disjoint supports.

To start with, observe that we can assume w.l.o.g. that each
$\Lambda_i(r)$ is bounded in $r$, otherwise \eq{tesi} would be
already proved. By this boundedness, we derive that, for $r$
large,
\neweq{controldenom}
\int_{S^{N-1}} (v_i^{(r)})^2\geq C>0.
\endeq
Indeed, assume not. This means that as $r\to\infty$ we have (up to
a subsequence) $\int_{S^{N-1}} (v_i^{(r)})^2\to 0$.
By the Holder inequality and since  $r$
is large, we infer that
$$\frac{1}{|\partial B_r|}\int_{\partial B_r} u_i\to 0.$$
Now we recall that $u_i$ is subharmonic, since it solves equation
\eq{k_nonsimm_intr}: hence, by the Mean Value Theorem and the
previous inequality, we have $u_i(0)=0$, a contradiction since
$u_i$ is strictly positive.

Let us now prove \eq{tesi} by showing that
$$\frac{2}{h}\sum_{i=1}^h
\sqrt{\Lambda_i(r)}>h'
\qquad \forall r\geq
r_0.$$
To this aim we argue by contradiction, and we assume
the existence of $r_n\to\infty$ such that
\neweq{fas}\frac{2}{h}\sum_{i=1}^h
\sqrt{\Lambda_i(r_n)}\leq   h'<\beta(h,N).\endeq
Let us define
$$
\tilde v_{i,n}=\frac{v_i^{(r_n)}}{\left(\int_{S^{N-1}} (v_i^{(r_n)})
^2\right)^{1/2}}.
$$
We have
\neweq{controlnumerat}
C\geq\Lambda_i(r_n)=\int_{S^{N-1}}\big[|\nabla\tilde
v_{i,n}|^2+r_n^2\big(\prod_{j\neq i}\|v_j^{(r_n)}\|_2\big) (\tilde
v_{i,n})^2\sum_{j\neq i}a_{ij}\tilde v_{j,n}\big].
\endeq
Since by definition $\int_{S^{N-1}}\tilde v_{i,n}^2=1$,  equation
\eq{controlnumerat} implies that $\int_{S^{N-1}}\nabla(\tilde
v_{i,n})^2$ is bounded. We infer the existence of a subsequence
$n_k\to\infty$ and $\bar{v}_i\not\equiv0$ such that $\tilde
v_{i,n_k}\rightharpoonup\bar{v}_i$, weakly in $H^1(S^{N-1})$, as
$k\to\infty$. This immediately gives
$$
\lim_{k\to\infty}\Lambda_i(r_{n_k})\geq \frac{\int_{S^{N-1}}|\nabla
\bar{v}_i|^2}{\int_{S^{N-1}}\bar{v}_i^2}=\l_1(\{\bar{v}_i>0\}).
$$
Taking into account \eq{controldenom}, we infer by
\eq{controlnumerat} that $ \tilde v_{i,n_k}\tilde v_{j,n_k}\to0$,
and therefore, by weak convergence, $\bar{v}_i\cdot
\bar{v}_j\equiv 0$ if $i\neq j$. Hence the supports of the
$\bar{v}_i$'s constitute a partition of $S^{N-1}$. Using this
information and summing up on $i$ the last inequality, we have
$$
\lim_{k\to\infty}\sum_{i=1}^h \sqrt{\Lambda_i(r_{n_k})}\geq \frac{h}{2}\beta(h,N).
$$
This provides a contradiction with \eq{fas}, concluding the proof.
\endproof

Let us conclude this section by pointing out a straightforward
consequence of the above monotonicity formula. Indeed, it induces some
growth restriction to the
solutions of \eq{k_nonsimm_intr} with positive components, as the following argument proves.
First, notice that Lemma \ref{monot.formula.N_intr}
with the choice $h=k$ and any $k'<k$ gives $\Phi(r)\geq\Phi(r')$
for all $r\geq r'$. We can assume w.l.o.g. that  $\Phi(r')=1$ so
that
\neweq{kk'}\prod_{i=1}^{k}\int_{B(0,r)}\left(|\nabla
u_i(x)|^2+u_i^2(x)\sum_{j\neq
i}a_{ij}u_j(x)\right)dx>r^{kk'}\hspace{1cm}r>r'\;.\endeq Let us
now go back to the differential equation for $u_i$: multiplying by
$u_i$ and integrating we have
$$
\int_{B(0,r)}\left(|\nabla u_i|^2+u_i^2\sum_{j\neq
i}a_{ij}u_j\right)=\int_{\partial B(0,r)} u_i\partial_n u_i\;.
$$
Let us now suppose that there exists $\a>0$ such that, for all $i$
$$
\partial_r u_i\leq C r^\a.
$$
Then, the r.h.s. is asymptotic to $r^{N-1}\cdot r^{\a+1}\cdot
r^\a$. Using this in \eq{kk'} we have, for $r$ large,
$$
r^{kk'}<\prod_{i=1}^{k}\int_{B(0,r)}\left(|\nabla
u_i|^2+u_i^2\sum_{j\neq i}a_{ij}u_j\right)\leq C r^{k(2\a+N)}.
$$
This provides $\a\geq(k'-N)/2$ for every $k'<k$. Hence we can
state the following
\begin{prop}\label{growth}
Let $U=(u_1,...,u_k)$ be  a solution of \eq{k_nonsimm_intr} on
$\bbbr^N$ with strictly positive components. Assume that there exists
$\a>0$ such that $|\nabla U|\leq Cr^\a$. Then  $\a\geq
(k-N)/2$.
\end{prop}
In particular, if $k\geq N$ then all the positive solutions
of \eq{k_nonsimm_intr} are unbounded at infinity, together with
their gradients.

\section{Liouville type Theorems}

In the spirit of Proposition \ref{growth}, let us focus our
attention on some nonexistence results of Liouville--type which
follow by application of the monotonicity formulae. We start by
proving that the system \eq{k_nonsimm_intr} does not admit
H\"older continuous solutions: this is a crucial step when
analyzing the rate of convergence of a class of
competition--diffusion systems, as the parameter of interspecific
competition tends to infinity (\cite{ctv4}).
\begin{prop}\label{liuv_perturb}
Let $k\geq 2$ and let $U$ be a solution of \eq{k_nonsimm_intr} on
$\bbbr^N$. Let $\a\in (0,1)$ such that
\neweq{hold}
\max_{i=1,\dots,k}\sup_{{x\in \bbbr^N}}
\frac{|u_i(x)|}{1+|x|^\a}<\infty.
\endeq
Then, $k-1$ components annihilate
and the last is a nonnegative constant.
\end{prop}

\proof by the strong maximum principle, every $u_i$ is either
identically zero or strictly positive. Let $h$ be the number of
the components not identically zero. If $h=1$, then the
proposition follows, without any growth assumption, by the
straight application of Liouville's theorem. Hence let $h\geq2$
and $u_1,\dots,u_h$ strictly positive. Let $B_r=B(0,r)$. By Lemma
\ref{monot.formula.N_intr} we know that
\neweq{so}\prod_{i=1}^h\frac{1}{r^{h'}}\int_{B_r} \left(|\nabla
u_i(x)|^2+u_i^2(x)\sum_{j\neq i}a_{ij}u_j(x)\right)dx\geq C>0
\endeq
when $h'<\beta(h,N)$ and $r$ is large enough.

On the other hand, let us consider a smooth, radial cut--off
function which is equal 1 in $B_r$ and vanishes outside $B_{2r}$.
Let us multiply the $i$--th differential equation by $\eta^2u_i$
and then integrate:
$$
\begin{array}{c}
\int_{B_{2r}}\eta^2|\nabla u_i|^2+\eta^2 u_i^2\sum_{j\neq i}
u_j\leq \int_{B_{2r}}|2\eta u\nabla\eta\nabla u|\\ \leq
\frac{1}{2} \int_{B_{2r}}\eta^2|\nabla u_i|^2 +2
\int_{B_{2r}}u_i^2|\nabla \eta|^2
\end{array}
$$
and hence
$$
\int_{B_{r}}|\nabla u_i|^2+u_i^2\sum_{j\neq i} u_j\leq
4\int_{B_{2r}\setminus B_r}u_i^2|\nabla \eta|^2\leq \frac{C}{r^2}
\int_{B_{2r}\setminus B_r}u_i^2.
$$
By \eq{hold}, when $\rho$ is sufficiently large, we have that
$u(x)\leq C' \rho^\a$ for all $x\in\partial B_{\rho}$: using this
fact in the above inequality and passing to polar coordinates we
obtain
$$
\int_{B_{r}}|\nabla u_i|^2\leq \frac{C}{r^2}\int_r^{2r}
\rho^{N-1+2\a}d\rho=C r^{N-2(1-\a)}
$$
for all indices $i$. Comparing with \eq{so} we have $r^{h'}\leq
Cr^{N-2(1-\a)}$ for $r$ large enough. But now, using Proposition
\ref{beta.cresce}, we can choose $h':=N-(1-\a)<N\leq\beta(h,N)$,
which provides a contradiction.
\endproof
Following the same line of the previous proof, but exploiting the
monotonicity formula Lemma \ref{monot.generale_intr}
(with $\bar{r}=\infty$), the subsequent result follows at once:
\begin{prop}\label{liuv_subarmonica}
Let $k\geq 2$ and let $U=(u_1,\dots,u_k)$ such that $u_i\cdot
u_j=0$ if $i\neq j$ and $-\Delta u_i\leq 0$ on $\bbbr^N$ for all
$i$. Let $\a\in (0,1)$ such that
$$\sup_{x\in \bbbr^N}
\frac{|u_i(x)|}{1+|x|^\a}<\infty$$ for all $i=1,\dots, k$. Then
each component $u_i$ is constant.
\end{prop}
It is worthwhile noticing that an analogous nonexistence result
holds for harmonic functions on the entire space and it can be
proved in a similar fashion by using the original monotonicity
formula Lemma \ref{clax_intr}.
\begin{prop}\label{liuv_armonica}
Let $u$ be an harmonic function on $\bbbr^N$. Let $\a\in (0,1)$
such that
$$\sup_{x\in \bbbr^N}
\frac{|u(x)|}{1+|x|^\a}<\infty.$$ Then $u$ is constant.
\end{prop}
In fact, the last assertion could be proved even easier, by a
simple test of the equation $-\Delta u=0$ with $u$.

\end{document}